# On Minkowski Inequalities Involving Fractional Calculus With General Analytic Kernels

Erdal GÜL[1], Ahmet Ocak AKDEMİR[2] and Abdüllatif YALÇIN [3,*]

[1]Yıldız Technical University, Faculty of Arts and Science
Department of Mathematics, Istanbul, Turkey

[2]Ağrı İbrahim Çeçen University, Faculty of Science and Letters, Department of Mathematics,
Ağrı, Turkey

[3] Yıldız Technical University, Graduate School of Applied and Natural Sciences,
Department of Mathematics, Istanbul, Turkey

∗Corresponding author: abdullatif.yalcin@std.yildiz.edu.tr

**Abstract:** There have been many proposed forms of fractional calculus, which can be grouped into a few broad classes of operators. By replacing the kernel of the power function with another kernel function, the traditional Riemann-Liouville formula and its generalisations are modified. Recent research has focused on Minkowski fractional inequalities and other inequalities for fractional integrals of certain types of functions. We also provide many concrete examples of applications to demonstrate the power of our key findings. It is no longer necessary to prove such statements separately for each model. This is because in a single study we have provided theorems that are valid for the entire general class of fractional operators. We can also study fractional integrals and derivatives with respect to functions.

**Key words:** Fractional integrals; Minkowski inequality; Integral inequalities; General analytic kernels.

1. Introduction

Through differential and integral calculations, calculus, one of the fundamental ideas of mathematics, is essential to many scientific and engineering disciplines. However, there are occasions when the conventional integer-order differentiation and integration processes are insufficient for the analysis of complex systems or the modelling of natural phenomena. For example, mathematicians have been fascinated by the $1/2$ order derivative for hundreds of years. At present, applications and research in mathematics are gaining a new perspective thanks to fractional calculus. However, there is still no conclusive answer to this question. In the 21st century, research in this area, which has its origins in the 17th century, has grown and deepened, and many more useful definitions have been developed. And each has its own advantages and disadvantages. One of the most useful fractional integrals is the Riemann-Liouville fractional integral operator (see [1-5]). By challenging mathematical thinking in new ways, fractional calculus has had a major impact on various fields of application. Fractional operators arise in two areas, pure mathematical theory and real-world modelling, and play an important role in disciplines such as engineering, physics, statistics, biology and economics. They also allow us to better understand and model complex phenomena (see [6-9]).

In recent years, researchers have introduced the Caputo model by modifying the differentiation and integration operations of the most popular model, the Riemann-liouville fractional integral operator (see [10]). Caputo is often preferred for modeling initial value problems, but the disadvantage of this model is that the analyticity is lost during differential integration (see [11]). Most of the definitions of fractional operators are based on replacing the kernel of the power function by a different kernel function. In particular, we refer to the Atangana Baleanu and Prabhakar definitions involving the Mitag-Lefler kernel function (which is nonsingular and interpolates between Riemann-Liouville and Caputo) (see [12-18,38]). (see [12-18]) These definitions are not equivalent to each other, each has its own specific properties and applications, with the Prabhakar formula attracting attention in recent years for its application in viscoelasticity and scotastic processes (see [19,20]). For example, to model systems with power law behavior, logarithmic behavior, exponential behavior and other more complex behaviors or for different initial conditions. An individual type of fractional calculus can be useful in different types of application domains. However, from a mathematical point of view, it is not efficient to prove the same properties and theorems over and over again. For each fractional operator, researchers try to generalize and prove results again and again.

It is a natural question to ask from a mathematical point of view whether it is possible to implement them only once in a general setting and then as special cases for each operator. A similar question has been raised by researchers working in applied sciences, such as engineers, in terms of real-world applications. Therefore, recently many results in the class of fractional integrals and derivative operators with analytic kernels have been proved in a general setting (see [21]).

In this work we try to stop the flow of similar papers by proving results about general classes. By eliminating the proofs of each class of fractional operators, we try to eliminate the problems of wasted time and class differences. We continue this paper with proofs involving the inverse Minkowski inequality and related results and the general class of fractional operators with analytic kernel.

The paper is organized as follows: In Section 2 we discuss notions of results and basic definitions of the newly introduced general analytic kernel fractional operators. We also present results on the inverse Minkowski inequality. Section 3 defends fundamental results such as the inverse Minkowski inequality involving integrals with general analytic kernels. Section 4, we show related variants. Section 5 is we offer many concrete examples as applications. Section 6 is devoted to conclusions.

2. Preliminaries

2.1 Fractional Calculus

We begin this chapter by describing the most popular point of the fractional integral, the Riemann-Liouville fractional integral operator, and fractional integral models based on replacing the kernel of a power function by a different kernel function.

**Definition 2.1 ([22, 23]).** Let $f \in \mathcal{L}^1[a,b]$ and $R(\alpha) > 0$ the right and left Riemann Liouville fractional integrals $^{RL}I_{a+}^\alpha f$ and $^{RL}I_{b-}^\alpha f$ of order $\alpha$ of $f(x)$ is are defined by:

$$^{RL}I_{a+}^\alpha f(x) = \frac{1}{\Gamma(\alpha)} \int_a^x f(\theta)(x-\theta)^{\alpha-1} d\theta, \quad x > a,$$

(2.1)

and

$$^{RL}I_{b-}^\alpha f(x) = \frac{1}{\Gamma(\alpha)} \int_x^b f(\theta)(\theta-x)^{\alpha-1} d\theta, \quad b > x.$$

(2.2)

Where $\Gamma(\alpha)$ is gamma function.

**Definition 2.2 ([5, 22, 23]).** Let $f \in \mathbb{C}^n[a,b]$ a function $n-1 \leq R(\alpha) < 0$. Then, the Riemann–Liouville fractional derivative with order $\alpha$ of the function $f$ with respect to $x$ and with constant of integration a is defined by:

$$^{RL}D_{a+}^\alpha f(x) = \frac{d^n}{dx^n}\left(^{RL}I_{a+}^{n-\alpha} f(x)\right),$$

(2.3)

and

$$^{RL}D^{\alpha}_{b^-}f(x) = \frac{d^n}{dx^n}\left(^{RL}I^{n-\alpha}_{b^-}f(x)\right).$$

By conventional, we have $^{RL}D^{-\alpha}_{a^+}f(x) = {}^{RL}I^{\alpha}_{a^+}$, we find that $^{RL}D^{\alpha}_{a^+}f(x)$ and $^{RL}I^{\alpha}_{a^+}$ are defined for all $\alpha \in \mathbb{C}$.

By differentiation and fractional integration of the kernel of the power function of the Riemann-Liouville fractional operator, i.e. another function, this gives rise to the Caputo and Atangana-Baleanu model. The Caputo-Fabrizio fractional operator, which is used in dynamical systems, physical phenomena, disease models, and many other fields, is a highly functional operator by definition, But has a deficiency in terms of not meeting the initial conditions in the special case $\alpha = 1$. The improvement to eliminate this deficiency has been provided by the new derivative operator developed by Atangana-Baleanu, which has versions in the sense of Caputo and Riemann (see [10,12,14,15]).

**Definition 2.2 ([24]).** Let $[a, b]$ be a real interval, $\alpha$ and $\beta$ be complex parameters with non-negative real parts, and $\mathbb{R}$ be a positive number satisfying $\mathbb{R} > (b-a)^{Re(\beta)}$ Let A be a complex function analytic on the disc $D(0, \mathbb{R})$ and defined on this disc by the locally uniformly convergent power series

$$A(x) = \sum_{n=0}^{\infty} a_n x^n,$$

(2.4)

where the coefficients $a_n = a_n(\alpha, \beta)$ are permitted to depend on $\alpha$ and $\beta$ if desired. We define the following fractional integral operatör $f \in \mathcal{L}[a, b]$ and $\alpha, \beta > 0$, the right and left Riemann Liouville fractional integral with analytic kernel function $^A I^{\alpha,\beta}_{a^+} f$ and $^A I^{\alpha,\beta}_{b^-} f$ of order $\alpha$ of $f(x)$ is are defined by:

$$^A I^{\alpha,\beta}_{a^+} f(x) = \int_a^x f(\theta)(x-\theta)^{\alpha-1} A\big((x-\theta)^\beta\big) d\theta, \quad x > a,$$

(2.5)

and

$$^A I^{\alpha,\beta}_{b^-} f(x) = \int_x^b f(\theta)(\theta-x)^{\alpha-1} A\big((t-\theta)^\beta\big) d\theta, \quad b > x.$$

(2.6)

The formula (2.5) and (2.6) is an extreme generalisation of the assortment of fractional models we considered above.

**Definition 2.3 ([24]).** For any analytic function A as in Definition 2.2, we define $A_\Gamma$ to be the transformed function

$$A_\Gamma(x) = \sum_{n=0}^{\infty} a_n \, \Gamma(\beta n + \alpha) x^n.$$

(2.7)

The relationship between the pair of functions A and $A_\Gamma$ is vital to the understanding of our generalised operators.

Alternatively, the generalised integral operator with analytical kernel (2.5) and (2.6) can be written as an infinite series of Riemann-Liouville fractional integrals. It is thus confirmed that the fractional calculus is part of the calculus by the following theorem expressed as follows (see [24]).

**Theorem 2.1 ([24]).** With all notation as in Definition 2.2, for any function $f \in \mathcal{L}^1[a, b]$, we have the following locally uniformly convergent series for ${}^A I_{a^+}^\alpha f(x)$ as a function on $[a, b]$:

$$^A I_{a^+}^{\alpha,\beta} f(x) = \sum_{n=0}^{\infty} a_n \, \Gamma(\beta n + \alpha) \, {}^{RL} I_{a^+}^{\alpha + n\beta} f(x)$$

Alternatively, this identity can be written more concisely in terms of the transformed function $A_\Gamma$ introduced in (2.7):

$$^A I_{a^+}^{\alpha,\beta} f(x) = A_\Gamma \left( {}^{RL} I_{a^+}^{\beta} f(x) \right) \left( {}^{RL} I_{a^+}^{\alpha} f(x) \right).$$

Similarly, the right fractional integral can also be defined with the analytic kernel function A and the parameters $\alpha, \beta$ of $f(x)$.

It has been shown that the Prabhakar kernel is general enough to include as special cases some other kernel functions of fractional calculus, including the AB kernel (see [13]). It has also been shown that all the fractional models mentioned above can be viewed as special cases of our generalised model with the new analytical kernel for appropriate functions $f$ and parameters $\alpha$, $\beta$, A (see [24]).

### 2.2 Minkowski inequality

The reverse Minkowski inequality and a result relevant to the inequality connected to the Riemann-Liouville fractional integral that corresponds to the following two theorems were proven by Dahmani (see [25]).

**Theorem 2.2 ([25]).** Let $\alpha > 0$, $p \geq 1$ and let there be two positive functions $f_1$ and $f_2$ on $[0, \infty)$ such that for all $x > a$, ${}^{RL}I_{a^+}^{\alpha} f_1^{\,p}(x) < \infty$, ${}^{RL}I_{a^+}^{\alpha} f_2^{\,p}(x) < \infty$. If $0 < \tau_1 \leq \frac{f_1(\theta)}{f_2(\theta)} \leq \tau_2$, $\theta \in [a, x]$, then we have:

$$\left({}^{RL}I_{a^+}^{\alpha} f_1^{\,p}(x)\right)^{\frac{1}{p}} + \left({}^{RL}I_{a^+}^{\alpha} f_2^{\,p}(x)\right)^{\frac{1}{p}} \leq \frac{1 + \tau_2(\tau_1 + 2)}{(\tau_1 + 1)(\tau_2 + 1)} \left({}^{RL}I_{a^+}^{\alpha} (f_1 + f_2)^p(x)\right)^{1/p}.$$

(2.7)

**Theorem 2.3 ([25]).** Let $\alpha > 0$, $p \geq 1$ and let there be two positive functions $f_1$ and $f_2$ on $[0, \infty)$ such that for all $x > a$, ${}^{RL}I_{a^+}^{\alpha} f_1^{\,p}(x) < \infty$, ${}^{RL}I_{a^+}^{\alpha} f_2^{\,p}(x) < \infty$. If $0 < \tau_1 \leq \frac{f_1(\theta)}{f_2(\theta)} \leq \tau_2$, $\theta \in [a, x]$, then we have:

$$\left({}^{RL}I_{a^+}^{\alpha} f_1^{\,p}(x)\right)^{\frac{2}{p}} + \left({}^{RL}I_{a^+}^{\alpha} f_2^{\,p}(x)\right)^{\frac{2}{p}}$$
$$\geq \left(\frac{(1 + \tau_2)(\tau_1 + 1)}{\tau_2} - 2\right) \left[{}^{RL}I_{a^+}^{\alpha} f_1^{\,p}(x)\right]^{1/p} \left[{}^{RL}I_{a^+}^{\alpha} f_2^{\,p}(x)\right]^{1/p}.$$

(2.8)

The classical version of the reverse Minkowski inequality was proved by Bougoffa in [28], the version with the Mittag-Leffler kernel by Andric et al. in [26], the version with the AB-fractional operator by Khan et al. in [27], and many other extensions of the reverse Minkowski inequality have been proved by various types of fractional calculus.(see [29-37]) We do not give all the proofs in detail here, but they are often functionally identical to the results of the inequalities in (2.7) and (2.8).

### 3 Minkowski inequalities for fractional integral with general analytic kernels

This section comprises our principal involvement of establishing the proof of the reverse Minkowski inequalities via general analytic kernels fractional integral operators defined in (2.5) and (2.6) and an associated theorem insinuated as the reverse Minkowski inequalities.

**Theorem 3.1.** Let $\alpha, \beta > 0$, $p \geq 1$ and let there be two positive functions $f_1$ and $f_2$ on $[0, \infty)$ such that for all $x > a$, ${}^A I_{a^+}^{\alpha,\beta} f_1^p(x) < \infty$, ${}^A I_{a^+}^{\alpha,\beta} f_2^p(x) < \infty$. If $0 < \tau_1 \leq \frac{f_1(\theta)}{f_2(\theta)} \leq \tau_2$, holds for $\tau_1, \tau_2 \in \mathbb{R}^+$ and $\theta \in [a, x]$, then we have:

$$\left({}^A I_{a^+}^{\alpha,\beta} f_1^p(x)\right)^{\frac{1}{p}} + \left({}^A I_{a^+}^{\alpha,\beta} f_2^p(x)\right)^{\frac{1}{p}} \leq \frac{1 + \tau_2(\tau_1 + 2)}{(\tau_1 + 1)(\tau_2 + 1)} \left({}^A I_{a^+}^{\alpha,\beta} (f_1 + f_2)^p(x)\right)^{\frac{1}{p}}.$$

(3.1)

*Proof:* Under the given condition $\frac{f_1(\theta)}{f_2(\theta)} \leq \tau_2$, $\theta \in [a, x]$, it can be write as

$$(\tau_2 + 1)^p f_1^p(\theta) \leq \tau_2^p (f_1 + f_2)^p(\theta).$$

(3.2)

Multiplying both sides of (3.2) with $\frac{(x - \theta)^{\alpha + n\beta - 1}}{\Gamma(\beta n + \alpha)}$, since the gamma function is positive $\mathbb{R}^+$, and $(x - \theta)^{\alpha + n\beta - 1}$ is also positive. Then integrating the resulting inequalities with respect to $\theta$ over $(a, x)$, we obtain

$$(\tau_2 + 1)^p \int_a^x (x - \theta)^{\alpha + n\beta - 1} f_1^p(\theta) d\theta \leq \tau_2^p \int_a^x (x - \theta)^{\alpha + n\beta - 1} (f_1 + f_2)^p(\theta) d\theta.$$

(3.3)

Consequently, we can write

$$(\tau_2 + 1)^p {}^{RL} I_{a^+}^{\alpha + n\beta} f_1^p(x) \leq \tau_2^p {}^{RL} I_{a^+}^{\alpha + n\beta} (f_1 + f_2)^p(x).$$

(3.4)

If both sides of (3.4) are multiplied by $\Gamma(\beta n + \alpha)$ and all $a_n$ real positives and then summed over all $n$:

$$(\tau_2 + 1)^p \sum_{n=0}^{\infty} a_n \Gamma(\beta n + \alpha) {}^{RL} I_{a^+}^{\alpha + n\beta} f_1^p(x) \leq \tau_2^p \sum_{n=0}^{\infty} a_n \Gamma(\beta n + \alpha) {}^{RL} I_{a^+}^{\alpha + n\beta} (f_1 + f_2)^p(x).$$

Which is equivalent to

$$^A I_{a^+}^{\alpha,\beta} f_1^p(x) \leq \frac{\tau_2^p}{(\tau_2 + 1)^p} \, ^A I_{a^+}^{\alpha,\beta} (f_1 + f_2)^p(x).$$

Hence, we can write

$$\left[ ^A I_{a^+}^{\alpha,\beta} f_1^p(x) \right]^{\frac{1}{p}} \leq \frac{\tau_2}{(\tau_2 + 1)} \left[ ^A I_{a^+}^{\alpha,\beta} (f_1 + f_2)^p(x) \right]^{\frac{1}{p}}.$$

(3.5)

In contrast, as $\tau_1 f_2(\theta) \leq f_1(\theta)$, it follows that

$$\left(1 + \frac{1}{\tau_1}\right)^p f_2^p(\theta) \leq \frac{1}{\tau_1^p} [f_1(\theta) + f_2(\theta)]^p.$$

(3.6)

Again, if we multiplying both sides of (3.6) with $\frac{(x-\theta)^{\alpha+n\beta-1}}{\Gamma(\beta n + \alpha)}$, and all $a_n$ real positives and then summed over all $n$ we obtain.

$$\left[ ^A I_{a^+}^{\alpha,\beta} f_2^p(x) \right]^{\frac{1}{p}} \leq \frac{1}{(\tau_1 + 1)} \left[ ^A I_{a^+}^{\alpha,\beta} (f_1 + f_2)^p(x) \right]^{\frac{1}{p}}.$$

(3.7)

Adding the inequalities (3.5) and (3.7) yields the desired inequality.

**Remark 3.2.** Applying Theorem 3.1 for $\alpha = 1$, $\beta = 0$ and for an arbitrary choice of function $A(1)$ we obtain Theorem 1.2 in [28].

**Remark 3.3.** In Theorem 3.1, if we choose , $\beta = 0$ and $A(x) = \frac{1}{\Gamma(\alpha)}$, we obtain Theorem 2.1 in [25].

**Remark 3.4.** Using Theorem 3.1 with $A(x) = \frac{1}{\rho^\alpha \Gamma(\alpha)} \exp\left(\frac{\rho-1}{\rho} x\right)$, we obtain Theorem 3.1 in [36].

**Remark 3.5.** In Theorem 3.1, if we choose $A(x) = E_{\beta,\alpha}^\rho(\omega x)$, we obtain Theorem 2.1 in [26].

Inequality (3.1) is referred to as the reverse Minkowski inequality for fractional integral with general analytic kernels.

**Theorem 3.2.** Let $\alpha, \beta > 0$, $p \geq 1$ and let there be two positive functions $f_1$ and $f_2$ on $[0, \infty)$ such that for all $x > a$, ${}^A I_{a^+}^{\alpha,\beta} f_1^p(x) < \infty$, ${}^A I_{a^+}^{\alpha,\beta} f_2^p(x) < \infty$. If $0 < \tau_1 \leq \frac{f_1(\theta)}{f_2(\theta)} \leq \tau_2$, holds for $\tau_1, \tau_2 \in \mathbb{R}^+$ and $\theta \in [a, x]$, then we have:

$$\left({}^A I_{a^+}^{\alpha,\beta} f_1^p(x)\right)^{\frac{2}{p}} + \left({}^A I_{a^+}^{\alpha,\beta} f_2^p(x)\right)^{\frac{2}{p}}$$
$$\geq \left(\frac{(1+\tau_2)(\tau_1+1)}{\tau_2} - 2\right) \left[{}^A I_{a^+}^{\alpha,\beta} f_1^p(x)\right]^{\frac{1}{p}} \left[{}^A I_{a^+}^{\alpha,\beta} f_1^p(x)\right]^{\frac{1}{p}}.$$

(3.8)

*Proof:* The product of inequalities (3.5) and (3.7) yields

$$\frac{(1+\tau_2)(\tau_1+1)}{\tau_2} \left[{}^A I_{a^+}^{\alpha,\beta} f_1^p(x)\right]^{\frac{1}{p}} \left[{}^A I_{a^+}^{\alpha,\beta} f_2^p(x)\right]^{\frac{1}{p}} \leq \left[{}^A I_{a^+}^{\alpha,\beta} (f_1+f_2)^p(x)\right]^{\frac{2}{p}}.$$

(3.9)

Now, utilizing the Minkowski inequality to the right hand side of (3.8), one obtains

$$\left({}^A I_{a^+}^{\alpha,\beta} (f_1+f_2)^p(x)\right)^{\frac{2}{p}} \leq \left(\left[{}^A I_{a^+}^{\alpha,\beta} f_1^p(x)\right]^{\frac{1}{p}} + \left[{}^A I_{a^+}^{\alpha,\beta} f_2^p(x)\right]^{\frac{1}{p}}\right)^2.$$

Then, we have

$$\left({}^A I_{a^+}^{\alpha,\beta} (f_1+f_2)^p(x)\right)^{\frac{2}{p}} \leq \left[{}^A I_{a^+}^{\alpha,\beta} f_1^p(x)\right]^{\frac{2}{p}} + \left[{}^A I_{a^+}^{\alpha,\beta} f_2^p(x)\right]^{\frac{2}{p}} + 2\left[{}^A I_{a^+}^{\alpha,\beta} f_1^p(x)\right]\left[{}^A I_{a^+}^{\alpha,\beta} f_2^p(x)\right].$$

(3.10)

Thus, from inequalities (3.9) and (3.10), we obtain the inequality (3.8).

**Remark 3.6.** Applying Theorem 3.1 for $\alpha = 1$, $\beta = 0$ and for an arbitrary choice of function $A(1)$ we obtain Theorem 2.2 in [28].

**Remark 3.7.** In Theorem 3.1, if we choose , $\beta = 0$ and $A(x) = \frac{1}{\Gamma(\alpha)}$, we obtain Theorem 2.3 in [25].

**Remark 3.8.** Using Theorem 3.1 with $A(x) = \frac{1}{\rho^\alpha \Gamma(\alpha)} exp\left(\frac{\rho-1}{\rho} x\right)$, we obtain Theorem 3.2 in [36].

**Remark 3.9.** In Theorem 3.1, if we choose $A(x) = E_{\beta,\alpha}^{\rho}(\omega x)$, we obtain Theorem 2.2 in [26].

## 4 Certain associated inequalities via fractional integral for general analytic kernels

This section is devoted to deriving certain related inequalities involving a generalized with analaytic kernels fractional integral operator.

**Theorem 4.1.** Let $\alpha, \beta > 0$, $p \geq 1$ and $\frac{1}{p} + \frac{1}{q} = 1$ let there be two positive functions $f_1$ and $f_2$ on $[0, \infty)$ such that for all $x > a$, $^A I_{a^+}^{\alpha,\beta} f_1^p(x) < \infty$, $^A I_{a^+}^{\alpha,\beta} f_2^p(x) < \infty$. If $0 < \tau_1 \leq \frac{f_1(\theta)}{f_2(\theta)} \leq \tau_2$, holds for $\tau_1, \tau_2 \in \mathbb{R}^+$ and $\theta \in [a, x]$, then we have:

$$\left[{}^A I_{a^+}^{\alpha,\beta} f_1(x)\right]^{\frac{1}{p}} \left[{}^A I_{a^+}^{\alpha,\beta} f_2(x)\right]^{\frac{1}{q}} \leq \left(\frac{\tau_2}{\tau_1}\right)^{\frac{1}{pq}} \left[{}^A I_{a^+}^{\alpha,\beta} f_1^{\frac{1}{p}}(x) f_2^{\frac{1}{q}}(x)\right].$$

(4.1)

*Proof:* Under the given condition $\frac{f_1(\theta)}{f_2(\theta)} \leq \tau_2$, $\theta \in [a, x]$, it can be write as

$$f_1(\theta) \leq \tau_2 f_2(\theta)$$

$$\tau_2^{-\frac{1}{q}} f_1^{\frac{1}{q}}(\theta) \leq f_2^{\frac{1}{q}}.$$

(4.2)

Taking the product of both sides of (4.2) by $f_1^{\frac{1}{p}}$, we can write as follows

$$\tau_2^{-\frac{1}{q}} f_1(\theta) \leq f_2^{\frac{1}{q}}(\theta) f_1^{\frac{1}{p}}(\theta).$$

(4.3)

Just like in the proof of Theorem 3.1, multiplying both sides of (4.3) with $\frac{(x-\theta)^{\alpha+n\beta-1}}{\Gamma(\beta n+\alpha)}$, and all $a_n$ real positives and then summed over all $n$ and then integrating the resulting inequalities with respect to $\theta$ over $(a, x)$, we obtain the following inequalities.

$$\tau_2^{-\frac{1}{q}} \sum_{n=0}^{\infty} a_n \Gamma(\beta n + \alpha) {}^{RL}I_{a^+}^{\alpha+n\beta} f_1(x) \leq \sum_{n=0}^{\infty} a_n \Gamma(\beta n + \alpha) {}^{RL}I_{a^+}^{\alpha+n\beta} f_2^{\frac{1}{q}}(x) f_1^{\frac{1}{p}}(x).$$

Consequently, we have

$$\tau_2^{-\frac{1}{pq}} \left[ {}^A I_{a^+}^{\alpha,\beta} f_1(x) \right]^{\frac{1}{p}} \leq \left[ {}^A I_{a^+}^{\alpha,\beta} f_1^{\frac{1}{p}}(x) f_2^{\frac{1}{q}}(x) \right]^{\frac{1}{p}}.$$

(4.4)

On the contary, as $\tau_1 \leq \frac{f_1(\theta)}{f_2(\theta)}$, we have

$$\tau_1^{\frac{1}{p}} f_2^{\frac{1}{p}}(\theta) \leq f_1^{\frac{1}{p}}(\theta).$$

(4.5)

Multiplying both sides of (4.5) by $f_2^{\frac{1}{q}}(\theta)$ and invoke the relation $\frac{1}{p} + \frac{1}{q} = 1$, we obtain

$$\tau_1^{\frac{1}{p}} f_2(\theta) \leq f_1^{\frac{1}{p}}(\theta) f_2^{\frac{1}{q}}(\theta)$$

(4.6)

Multiplying both sides of (4.6) with $\frac{(x-\theta)^{\alpha+n\beta-1}}{\Gamma(\beta n+\alpha)}$, and all $a_n$ real positives and then summed over all $n$ and then integrating the resulting inequalities with respect to $\theta$ over $(a, x)$, we obtain

$$\tau_1^{\frac{1}{pq}} \left[ {}^A I_{a^+}^{\alpha,\beta} f_2(x) \right]^{\frac{1}{q}} \leq \left[ {}^A I_{a^+}^{\alpha,\beta} f_1^{\frac{1}{p}}(x) f_2^{\frac{1}{q}}(x) \right]^{\frac{1}{q}}.$$

(4.7)

Multiplying (4.4) and (4.7), the required inequality (4.1) can be concluded.

**Theorem 4.2.** Let $\alpha, \beta > 0$, $p \geq 1$ and $\frac{1}{p} + \frac{1}{q} = 1$ let there be two positive functions $f_1$ and $f_2$ on $[0, \infty)$ such that for all $x > a$, ${}^A I_{a^+}^{\alpha,\beta} f_1^p(x) < \infty$, ${}^A I_{a^+}^{\alpha,\beta} f_2^p(x) < \infty$. If $0 < \tau_1 \leq \frac{f_1(\theta)}{f_2(\theta)} \leq \tau_2$, holds for $\tau_1, \tau_2 \in \mathbb{R}^+$ and $\theta \in [a, x]$, then we have:

$$\left( {}^A I_{a^+}^{\alpha,\beta} f_1 f_2 \right)(x) \leq C_1 \left( {}^A I_{a^+}^{\alpha,\beta} (f_1 + f_2)^p \right)(x) + C_2 \left( {}^A I_{a^+}^{\alpha,\beta} (f_1 + f_2)^q \right)(x)$$

With $C_1 = \frac{2^{p-1} \tau_2^p}{p(\tau_2+1)^p}$ and $C_2 = \frac{2^{q-1}}{q(\tau_1+1)^q}$.

(4.8)

*Proof:* Using the hypothesis, we have

$$(\tau_2 + 1)^p f_1^p(\theta) \leq \tau_2^p (f_1 + f_2)^p(\theta).$$

(4.9)

Multiplying both sides of (4.9) with $\frac{(x-\theta)^{\alpha+n\beta-1}}{\Gamma(\beta n+\alpha)}$, and all $a_n$ real positives and then summed over all $n$ and then integrating the resulting inequalities with respect to $\theta$ over $(a, x)$, we obtain the following inequalities.

$$(\tau_2 + 1)^p \sum_{n=0}^{\infty} a_n \Gamma(\beta n + \alpha)^{RL}I_{a^+}^{\alpha+n\beta} f_1^p(x) \leq \tau_2^p \sum_{n=0}^{\infty} a_n \Gamma(\beta n + \alpha)^{RL}I_{a^+}^{\alpha+n\beta}(f_1 + f_2)^p(x).$$

Consequently, we have

$$^AI_{a^+}^{\alpha,\beta} f_1^p(x) \leq \frac{\tau_2^p}{(\tau_2 + 1)^p} \, ^AI_{a^+}^{\alpha,\beta}(f_1 + f_2)^p(x).$$

(4.10)

In contrast, using $0 < \tau_1 \leq \frac{f_1(\theta)}{f_2(\theta)}, \theta \in [a, x]$ holds we get

$$(\tau_1 + 1)^q f_2^q(\theta) \leq (f_1 + f_2)^q(\theta).$$

(4.11)

Similarly, multiplying both sides of (4.11) with $\frac{(x-\theta)^{\alpha+n\beta-1}}{\Gamma(\beta n+\alpha)}$, and all $a_n$ real positives and then summed over all $n$ and then integrating the resulting inequalities with respect to $\theta$ over $(a, x)$, we can write

$$^AI_{a^+}^{\alpha,\beta} f_2^q(x) \leq \frac{1}{(\tau_1 + 1)^q} \left( ^AI_{a^+}^{\alpha,\beta}(f_1 + f_2)^q(x) \right)$$

(4.12)

Now, taking into account Young's inequality,

$$f_1(\theta)f_2(\theta) \leq \frac{1}{p} f_1^p(x) + \frac{1}{q} f_2^q(x),$$

(4.13)

again, multiplying both sides of (4.13) with $\frac{(x-\theta)^{\alpha+n\beta-1}}{\Gamma(\beta n+\alpha)}$, and all $a_n$ real positives and then summed over all $n$ and then integrating the resulting inequalities with respect to $\theta$ over $(a, x)$, we obtain

$$\left( {}^{A}I_{a^+}^{\alpha,\beta} f_1(x)f_2(x) \right) \le \frac{1}{p}\left( {}^{A}I_{a^+}^{\alpha,\beta} f_1(x) \right) + \frac{1}{q}\left( {}^{A}I_{a^+}^{\alpha,\beta} f_2(x) \right).$$

(4.14)

Applying inequalities (4.10) and (4.12) in (4.14), we obtain

$$\left( {}^{A}I_{a^+}^{\alpha,\beta} f_1(x)f_2(x) \right) \le \frac{\tau_2^{\,p}}{p(\tau_2+1)^p}\left( {}^{A}I_{a^+}^{\alpha,\beta}(f_1+f_2)^p(x) \right) + \frac{1}{q(\tau_1+1)^q}\left( {}^{A}I_{a^+}^{\alpha,\beta}(f_1+f_2)^q(x) \right).$$

(4.15)

Now, using the inequality $(\hbar_1 + \hbar_2)^\ell \le 2^{\ell-1}(\hbar_1^{\,\ell} + \hbar_2^{\,\ell})$, $\ell \ge 1$ with $\hbar_1, \hbar_2 > 0$ in (4.15)

We have

$$\left( {}^{A}I_{a^+}^{\alpha,\beta} f_1 f_2 \right)(x) \le \frac{2^{p-1}\tau_2^{\,p}}{p(\tau_2+1)^p}\left( {}^{A}I_{a^+}^{\alpha,\beta}(f_1+f_2)^p \right)(x) + \frac{2^{q-1}}{q(\tau_1+1)^q}\left( {}^{A}I_{a^+}^{\alpha,\beta}(f_1+f_2)^q \right)(x).$$

This is the required result.

**Theorem 4.3.** Let $\alpha, \beta > 0, p \ge 1$. Let $f_1, f_2 \in L[a,x]$ be two positive functions on $[0,\infty)$ such that for all $x > a$, ${}^{A}I_{a^+}^{\alpha,\beta} f_1^{\,p}(x) < \infty$, ${}^{A}I_{a^+}^{\alpha,\beta} f_2^{\,p}(x) < \infty$. If $0 < \tau_1 \le \frac{f_1(\theta)}{f_2(\theta)} \le \tau_2$, holds for $\tau_1, \tau_2 \in \mathbb{R}^+$ and $\theta \in [a,x]$, then:

$$\frac{\tau_2+1}{\tau_2-\varphi}\left( {}^{A}I_{a^+}^{\alpha,\beta}[f_1-\varphi f_2]^p \right)^{\frac{1}{p}}(x) \le \left( {}^{A}I_{a^+}^{\alpha,\beta} f_1^{\,p} \right)^{\frac{1}{p}}(x) + \left( {}^{A}I_{a^+}^{\alpha,\beta} f_2^{\,p} \right)^{\frac{1}{p}}(x)$$

$$\le \frac{\tau_1+1}{\tau_1-\varphi}\left( {}^{A}I_{a^+}^{\alpha,\beta}[f_1-\varphi f_2]^p \right)^{\frac{1}{p}}(x).$$

(4.16)

*Proof:* Using the hypothesis, $0 < \varphi < \tau_1 \le \frac{f_1(\theta)}{f_2(\theta)} \le \tau_2$, we have

$$\varphi\tau_1 \le \varphi\tau_2 \Rightarrow \tau_1\varphi + \tau_1 \le \varphi\tau_1 + \tau_2 \le \tau_2\varphi + \tau_2$$

$$\Rightarrow (\tau_2+1)(\tau_1-\varphi) \le (\tau_1+1)(\tau_2-\varphi).$$

It follows that

$$\frac{(\tau_2+1)}{(\tau_2-\varphi)} \le \frac{(\tau_1+1)}{(\tau_1-\varphi)}.$$

Furthermore, we obtain

$$\tau_1 - \varphi \leq \frac{f_1(\theta) - \varphi f_2(\theta)}{f_2(\theta)} \leq \tau_2 - \varphi$$

Hence, we obtain

$$\frac{[f_1(\theta) - \varphi f_2(\theta)]^p}{(\tau_2 - \varphi)^p} \leq f_2^{\,p}(\theta) \leq \frac{[f_1(\theta) - \varphi f_2(\theta)]^p}{(\tau_1 - \varphi)^p}.$$

(4.17)

Moreover, we have

$$\frac{1}{\tau_2} \leq \frac{f_2(\theta)}{f_1(\theta)} \leq \frac{1}{\tau_1} \Longrightarrow \frac{\tau_1 - \varphi}{\varphi \tau_1} \leq \frac{f_1(\theta) - \varphi f_2(\theta)}{\varphi f_1(\theta)} \leq \frac{\tau_2 - \varphi}{\varphi \tau_2}$$

Which implies that

$$\left(\frac{\tau_2}{\tau_2 - \varphi}\right)^p \leq [f_1(\theta) - \varphi f_2(\theta)]^p \leq f_1^{\,p}(\theta) \leq \left(\frac{\tau_1}{\tau_1 - \varphi}\right)^p [f_1(\theta) - \varphi f_2(\theta)]^p.$$

(4.18)

Multiplying both sides of (4.17) with $\frac{(x-\theta)^{\alpha+n\beta-1}}{\Gamma(\beta n+\alpha)}$, and all $a_n$ real positives and then summed over all $n$ and then integrating the resulting inequalities with respect to $\theta$ over $(a, x)$, we obtain the following inequalities.

$$\left(\frac{1}{\tau_2 - \varphi}\right)^p \sum_{n=0}^{\infty} a_n \Gamma(\beta n + \alpha)^{RL} I_{a^+}^{\alpha+n\beta}[f_1(x) - \varphi f_2(x)]^p \leq \sum_{n=0}^{\infty} a_n \Gamma(\beta n + \alpha)^{RL} I_{a^+}^{\alpha+n\beta} f_2^{\,p}(x)$$
$$\leq \left(\frac{1}{\tau_1 - \varphi}\right)^p \sum_{n=0}^{\infty} a_n \Gamma(\beta n + \alpha)^{RL} I_{a^+}^{\alpha+n\beta} (f_1 + \varphi f_2)^p(x).$$

Consequently, we have

$$\frac{1}{\tau_2 - \varphi} \left(^A I_{a^+}^{\alpha,\beta}[f_1 - \varphi f_2]^p\right)^{\frac{1}{p}}(x) \leq \left(^A I_{a^+}^{\alpha,\beta} f_2^{\,p}\right)^{\frac{1}{p}}(x) \leq \frac{1}{\tau_1 - \varphi} \left(^A I_{a^+}^{\alpha,\beta}[f_1 - \varphi f_2]^p\right)^{\frac{1}{p}}(x).$$

(4.19)

Adopting the same technique with (4.18), one obtains

$$\frac{\tau_2}{\tau_2 - \varphi} \left(^A I_{a^+}^{\alpha,\beta}[f_1 - \varphi f_2]^p\right)^{\frac{1}{p}}(x) \leq \left(^A I_{a^+}^{\alpha,\beta} f_1^{\,p}\right)^{\frac{1}{p}}(x) \leq \frac{\tau_1}{\tau_1 - \varphi} \left(^A I_{a^+}^{\alpha,\beta}[f_1 - \varphi f_2]^p\right)^{\frac{1}{p}}(x).$$

(4.20)

Finally, by adding inequalities (4.19) and (4.20), we complete the proof of Theorem 4.3.

**Theorem 4.4.** Let $\alpha, \beta > 0, p \geq 1$. Let $f_1, f_2 \in L[a, x]$ be two positive functions on $[0, \infty)$ such that for all $x > a$, ${}^A I_{a^+}^{\alpha,\beta} f_1^p(x) < \infty$, ${}^A I_{a^+}^{\alpha,\beta} f_2^p(x) < \infty$. If $0 \leq m \leq f_1(\theta) \leq M$ and $0 \leq n \leq f_2(\theta) \leq N$, $\theta \in [a, x]$, then:

$$\left({}^A I_{a^+}^{\alpha,\beta} f_1^p\right)^{\frac{1}{p}}(x) + \left({}^A I_{a^+}^{\alpha,\beta} f_2^p\right)^{\frac{1}{p}}(x) \leq C_3 \left({}^A I_{a^+}^{\alpha,\beta} [f_1 + f_2]^p\right)^{\frac{1}{p}}(x).$$

(4.21)

With $C_3 = \dfrac{M(m+N) + N(n+M)}{(m+N)(n+M)}$.

*Proof:* Under the given hypothesis, it follows that

$$\frac{1}{N} \leq \frac{1}{f_2(\theta)} \leq \frac{1}{n}.$$

(4.22)

Conducting the product between (4.22) and $0 < m \leq f_1(\theta) \leq M$, we have

$$\frac{m}{N} \leq \frac{f_1(\theta)}{f_2(\theta)} \leq \frac{M}{n}.$$

(4.23)

From (4.23), we get

$$f_2^p(\theta) \leq \left(\frac{N}{m+N}\right)^p [f_1(\theta) + f_2(\theta)]^p$$

(4.24)

and

$$f_1^p(\theta) \leq \left(\frac{M}{n+M}\right)^p [f_1(\theta) + f_2(\theta)]^p.$$

(4.25)

Multiplying both sides of (4.24) and (4.24) with $\frac{(x-\theta)^{\alpha+n\beta-1}}{\Gamma(\beta n+\alpha)}$, and all $a_n$ real positives and then summed over all $n$ and then integrating the resulting inequalities with respect to $\theta$ over $(a, x)$, we obtain the following inequalities.

$$\left({}^{A}I_{a^+}^{\alpha,\beta} f_2^p\right)^{\frac{1}{p}}(x) \leq \frac{\mathcal{N}}{m+\mathcal{N}} \left({}^{A}I_{a^+}^{\alpha,\beta} [f_1+f_2]^p\right)^{\frac{1}{p}}(x)$$

(4.26)

and

$$\left({}^{A}I_{a^+}^{\alpha,\beta} f_1^p\right)^{\frac{1}{p}}(x) \leq \frac{\mathcal{M}}{n+\mathcal{M}} \left({}^{A}I_{a^+}^{\alpha,\beta} [f_1+f_2]^p\right)^{\frac{1}{p}}(x).$$

(4.27)

By adding inequalities (4.26) and (4.27), we attain the inequality (4.21).

**Theorem 4.5.** Let $\alpha, \beta > 0, p \geq 1$. Let $f_1, f_2 \in L[a,x]$ be two positive functions on $[0,\infty)$ such that for all $x > a$, ${}^{A}I_{a^+}^{\alpha,\beta} f_1^p(x) < \infty$, ${}^{A}I_{a^+}^{\alpha,\beta} f_2^p(x) < \infty$. If $0 < \tau_1 \leq \frac{f_1(\theta)}{f_2(\theta)} \leq \tau_2$, holds for $\tau_1, \tau_2 \in \mathbb{R}^+$ and $\theta \in [a,x]$, then:

$$\frac{1}{\tau_2}\left({}^{A}I_{a^+}^{\alpha,\beta}(f_1 f_2)(x)\right) \leq \frac{1}{(\tau_1+1)(\tau_2+1)}\left({}^{A}I_{a^+}^{\alpha,\beta}(f_1+f_2)^2\right)(x) \leq \frac{1}{\tau_1}\left({}^{A}I_{a^+}^{\alpha,\beta}(f_1 f_2)(x)\right).$$

(4.28)

*Proof:* Under the given suppositions $0 < \tau_1 \leq \frac{f_1(\theta)}{f_2(\theta)} \leq \tau_2$, it follows that

$$f_2(\theta)(\tau_1+1) \leq f_1(\theta) + f_2(\theta) \leq f_2(\theta)(\tau_2+1).$$

(4.29)

Also it follows that $\frac{1}{\tau_2} \leq \frac{f_2(\theta)}{f_1(\theta)} \leq \frac{1}{\tau_1}$, which yields

$$f_1(\theta)\left(\frac{\tau_2+1}{\tau_2}\right) \leq f_1(\theta) + f_2(\theta) \leq f_1(\theta)\left(\frac{\tau_1+1}{\tau_1}\right).$$

(4.30)

The product of (4.29) and (4.30) gives

$$\frac{f_1(\theta)f_2(\theta)}{\tau_2} \leq \frac{(f_1(\theta)+f_2(\theta))^2}{(\tau_1+1)(\tau_2+1)} \leq \frac{f_1(\theta)f_2(\theta)}{\tau_1}.$$

(4.31)

Now, multiplying both sides of (4.31) by $\frac{(x-\theta)^{\alpha+n\beta-1}}{\Gamma(\beta n+\alpha)}$, and all $a_n$ real positives and then summed over all $n$ and then integrating the resulting inequalities with respect to $\theta$ over $(a, x)$, we obtain the following inequalities.

$$\frac{1}{\tau_2}\sum_{n=0}^{\infty} a_n \Gamma(\beta n+\alpha)^{RL}I_{a^+}^{\alpha+n\beta}(f_1 f_2)(x)$$

$$\leq \frac{1}{(\tau_1+1)(\tau_2+1)}\sum_{n=0}^{\infty} a_n \Gamma(\beta n+\alpha)^{RL}I_{a^+}^{\alpha+n\beta}(f_1+f_2)^2(x)$$

$$\leq \frac{1}{\tau_1}\sum_{n=0}^{\infty} a_n \Gamma(\beta n+\alpha)^{RL}I_{a^+}^{\alpha+n\beta}(f_1 f_2)(x).$$

One observes that

$$\frac{1}{\tau_2}\left(^A I_{a^+}^{\alpha,\beta}(f_1 f_2)(x)\right) \leq \frac{1}{(\tau_1+1)(\tau_2+1)}\left(^A I_{a^+}^{\alpha,\beta}(f_1+f_2)^2\right)(x) \leq \frac{1}{\tau_1}\left(^A I_{a^+}^{\alpha,\beta}(f_1 f_2)(x)\right).$$

which is the desired result.

**Theorem 4.6.** Let $\alpha, \beta > 0, p \geq 1$. Let $f_1, f_2 \in L[a, x]$ be two positive functions on $[0, \infty)$ such that for all $x > a$, $^A I_{a^+}^{\alpha,\beta} f_1^p(x) < \infty$, $^A I_{a^+}^{\alpha,\beta} f_2^p(x) < \infty$. If $0 < \tau_1 \leq \frac{f_1(\theta)}{f_2(\theta)} \leq \tau_2$, holds for $\tau_1, \tau_2 \in \mathbb{R}^+$ and $\theta \in [a, x]$, then:

$$\left(^A I_{a^+}^{\alpha,\beta} f_1^p(x)\right)^{\frac{1}{p}} + \left(^A I_{a^+}^{\alpha,\beta} f_2^p(x)\right)^{\frac{1}{p}} \leq 2\left(^A I_{a^+}^{\alpha,\beta} \Upsilon^p(f_1, f_2)\right)^{\frac{1}{p}}(x),$$

where $\Upsilon(f_1(\theta), f_2(\theta)) = \max\left\{\tau_2\left[\left(\frac{\tau_2}{\tau_1}+1\right)f_1(\theta) - \tau_2 f_2(\theta)\right], \frac{(\tau_2+\tau_1)f_2(\theta)-f_1(\theta)}{\tau_1}\right\}.$

(4.32)

*Proof:* Under the given condition $0 < \tau_1 \leq \frac{f_1(\theta)}{f_2(\theta)} \leq \tau_2$, $\theta \in [a, x]$, it can be write as

$$0 < \tau_1 \leq \tau_2 + \tau_1 - \frac{f_1(\theta)}{f_2(\theta)}$$

(4.33)

and

$$\tau_2 + \tau_1 - \frac{f_1(\theta)}{f_2(\theta)} \leq \tau_2.$$

(4.34)

Hence, using (4.33) and (4.34) we get

$$f_2(\theta) < \frac{(\tau_2 + \tau_1)f_2(\theta) - f_1(\theta)}{\tau_1} \leq \Upsilon(f_1(\theta), f_2(\theta)),$$

(4.35)

where $\Upsilon(f_1(\theta), f_2(\theta)) = \max\left\{\tau_2\left[\left(\frac{\tau_2}{\tau_1} + 1\right)f_1(\theta) - \tau_2 f_2(\theta)\right], \frac{(\tau_2+\tau_1)f_2(\theta)-f_1(\theta)}{\tau_1}\right\}$. Also from the given supposition $0 < \frac{1}{\tau_2} \leq \frac{f_2(\theta)}{f_1(\theta)} \leq \frac{1}{\tau_1}$. In this way, we have

$$\frac{1}{\tau_2} \leq \frac{1}{\tau_2} + \frac{1}{\tau_1} - \frac{f_2(\theta)}{f_1(\theta)},$$

(4.36)

and

$$\frac{1}{\tau_2} + \frac{1}{\tau_1} - \frac{f_2(\theta)}{f_1(\theta)} \leq \frac{1}{\tau_1}.$$

(4.36)

From (4.36) and (4.36), we obtain

$$\frac{1}{\tau_2} \leq \frac{\left(\frac{1}{\tau_2} + \frac{1}{\tau_1}\right)f_1(\theta) - f_2(\theta)}{f_1(\theta)} \leq \frac{1}{\tau_1},$$

(4.36)

implying

$$f_1(\theta) \leq \tau_2\left(\frac{1}{\tau_2} + \frac{1}{\tau_1}\right)f_1(\theta) - \tau_2 f_2(\theta)$$

$$= \left(\frac{\tau_2}{\tau_1} + 1\right)f_1(\theta) - \tau_2 f_2(\theta)$$

$$\leq \tau_2\left[\left(\frac{\tau_2}{\tau_1}+1\right)f_1(\theta)-\tau_2 f_2(\theta)\right]$$

$$\leq \Upsilon\bigl(f_1(\theta),f_2(\theta)\bigr)$$

(4.37)

From (4.35) and (4.37), we have

$$f_1{}^p(\theta)\leq \Upsilon^p\bigl(f_1(\theta),f_2(\theta)\bigr),$$

(4.38)

and

$$f_2{}^p(\theta)\leq \Upsilon^p\bigl(f_1(\theta),f_2(\theta)\bigr).$$

(4.39)

Now, multiplying both sides of (4.38) by $\frac{(x-\theta)^{\alpha+n\beta-1}}{\Gamma(\beta n+\alpha)}$, and all $a_n$ real positives and then summed over all $n$ and then integrating the resulting inequalities with respect to $\theta$ over $(a,x)$, we obtain.

$$\sum_{n=0}^{\infty} a_n\,\Gamma(\beta n+\alpha)\,{}^{RL}I_{a^+}^{\alpha+n\beta}f_1{}^p(x)\leq \sum_{n=0}^{\infty} a_n\,\Gamma(\beta n+\alpha)\,{}^{RL}I_{a^+}^{\alpha+n\beta}\Upsilon^p\bigl(f_1(x),f_2(x)\bigr)$$

Accordingly,

$$\left({}^{A}I_{a^+}^{\alpha,\beta}f_1{}^p(x)\right)^{\frac{1}{p}}\leq \left({}^{A}I_{a^+}^{\alpha,\beta}\Upsilon^p(f_1,f_2)\right)^{\frac{1}{p}}(x).$$

Adopting the same technique for (4.39), we have

$$\left({}^{A}I_{a^+}^{\alpha,\beta}f_2{}^p(x)\right)^{\frac{1}{p}}\leq \left({}^{A}I_{a^+}^{\alpha,\beta}\Upsilon^p(f_1,f_2)\right)^{\frac{1}{p}}(x).$$

(4.40)

Hence, adding them together, we complete the proof of Theorem 4.6.

## 5. Examples

Each of the following examples illustrates the findings of this paper in this section.

**Example 5.1.** Assume that $\alpha, \beta, \hbar > 0, p \geq 1$. Let $f_1, f_2 \in L[a, x]$ be two positive functions on $[0, \infty)$ such that for all $x > a \geq 1$ and $\theta \in [a, x]$, then we have:

$$\left[\left(^A I_{a^+}^{\alpha,\beta} \theta^p\right)(x)\right]^{\frac{1}{p}} + \left[\left(^A I_{a^+}^{\alpha,\beta} (\hbar + \theta)^p\right)(x)\right]^{\frac{1}{p}} \leq \frac{3\hbar + 4}{2(\hbar + 2)} \left[\left(^A I_{a^+}^{\alpha,\beta} (\hbar + 2\theta)^p\right)(x)\right]^{\frac{1}{p}}.$$

*Proof.* By selectingTaking $f_1(\theta) = \theta + \hbar$ and $f_2(\theta) = \hbar$, we have respectively $\tau_1 = 1$ and $\tau_2 = \hbar + 1$. Applying Theorem 3.1, the desired outcome.

**Example 5.2.** For any $\alpha, \beta, \hbar > 0, p \geq 1$. Let $f_1, f_2 \in L[a, x]$ be two positive functions on $[0, \infty)$ such that for all $x > a \geq 1$ and $\theta \in [a, x]$, then we have:

$$\left[\left(^A I_{a^+}^{\alpha,\beta} \theta^p\right)(x)\right]^{\frac{2}{p}} + \left[\left(^A I_{a^+}^{\alpha,\beta} (\hbar + \theta)^p\right)(x)\right]^{\frac{2}{p}}$$

$$\geq \left(\frac{2}{\hbar + 1}\right) \left[\left(^A I_{a^+}^{\alpha,\beta} \theta^p\right)(x)\right]^{\frac{1}{p}} \left[\left(^A I_{a^+}^{\alpha,\beta} (\hbar + \theta)^p\right)(x)\right]^{\frac{1}{p}}.$$

*Proof.* By putting $f_1(\theta) = \theta + \hbar$ and $f_2(\theta) = \hbar$, we have respectively $\tau_1 = 1$ and $\tau_2 = \hbar + 1$ Applying Theorem 3.2, the desired outcome.

**Example 5.3.** Suppose that $\alpha, \beta, \hbar > 0$, $p \geq 1$ and $\frac{1}{p} + \frac{1}{q} = 1$. Let $f_1, f_2 \in L[a, x]$ be two positive functions on $[0, \infty)$ such that for all $x > a \geq 1$ and $\theta \in [a, x]$, the following inequality holds:

$$\left[\left(^A I_{a^+}^{\alpha,\beta} (\theta + \hbar)\right)(x)\right]^{\frac{1}{p}} \left[\left(^A I_{a^+}^{\alpha,\beta} \theta\right)(x)\right]^{\frac{1}{q}} \leq (\hbar + 1)^{\frac{1}{pq}} \left(^A I_{a^+}^{\alpha,\beta} (\theta + \hbar)^{\frac{1}{p}} \theta^{\frac{1}{q}}\right)(x).$$

*Proof.* Choosing $f_1(\theta) = \theta + \hbar$ and $f_2(\theta) = \hbar$, we have respectively $\tau_1 = 1$ and $\tau_2 = \hbar + 1$ Applying Theorem 4.1, the desired outcome.

**Example 5.4.** Assume that, $\alpha, \beta, \hbar > 0, p \geq 1$ and $\frac{1}{p} + \frac{1}{q} = 1$. Let $f_1, f_2 \in L[a, x]$ be two positive functions on $[0, \infty)$ such that for all $x > a \geq 1$ and $\theta \in [a, x]$, the following inequality holds:

$$\left({}^A I_{a^+}^{\alpha,\beta}\theta(\theta+k)\right)(x)$$
$$\leq \frac{2^{p-1}(1+k)^p}{p(2+k)^p}\left({}^A I_{a^+}^{\alpha,\beta}(\theta^p+(\theta+k)^p)\right)(x)$$
$$+\frac{1}{2q}\left({}^A I_{a^+}^{\alpha,\beta}(\theta^p+(\theta+k)^p)\right)(x)$$

*Proof.* By setting $f_1(\theta) = \theta + k$ and $f_2(\theta) = k$, we have respectively $\tau_1 = 1$ and $\tau_2 = k+1$

Applying Theorem 4.2, we get the desired result.

**Example 5.5.** Assume that, $\alpha, \beta, k > 0, p \geq 1$. Let $f_1, f_2 \in L[a,x]$ be two positive functions on $[0, \infty)$ such that for all $x > a \geq 1$ and $\theta \in [a,x]$, the following inequality holds:

$$\frac{k+2}{k+1-\varphi}\left({}^A I_{a^+}^{\alpha,\beta}[\theta(1-\varphi)+k]^p\right)^{\frac{1}{p}}(x) \leq \left({}^A I_{a^+}^{\alpha,\beta}\theta^p\right)^{\frac{1}{p}}(x) + \left({}^A I_{a^+}^{\alpha,\beta}(\theta+k)^p\right)^{\frac{1}{p}}(x)$$
$$\leq \frac{2}{1-\varphi}\left({}^A I_{a^+}^{\alpha,\beta}[\theta(1-\varphi)+k]^p\right)^{\frac{1}{p}}(x).$$

*Proof.* By putting $f_1(\theta) = \theta + k$ and $f_2(\theta) = k$, we have respectively $\tau_1 = 1$ and $\tau_2 = k+1$

Applying Theorem 4.3, we get the desired result.

**Example 5.6.** Suppose that $\alpha, \beta, k > 0, p \geq 1$. Let $f_1, f_2 \in L[a,x]$ be two positive functions on $[0, \infty)$ such that for all $x > a \geq 1$ and $\theta \in [a,x]$, the following inequality holds:

$$\left({}^A I_{a^+}^{\alpha,\beta}\sin^{2p}\theta\right)^{\frac{1}{p}}(x) + \left({}^A I_{a^+}^{\alpha,\beta}\cos^{2p}\theta\right)^{\frac{1}{p}}(x) \leq 2\left({}^A I_{a^+}^{\alpha,\beta}1\right)^{\frac{1}{p}}(x).$$

*Proof.* Taking $f_1(\theta) = \sin^2\theta$ and $f_1(\theta) = \cos^2\theta$, we obtain respectively, $m = n = 0$ and $M = N = 1$. Applying Theorem 4.4, we get the desired result.

**Example 5.7.** For any, $\alpha, \beta, k > 0, p \geq 1$. Let $f_1, f_2 \in L[a,x]$ be two positive functions on $[0, \infty)$ such that for all $x > a \geq 1$ and $\theta \in [a,x]$, the following inequality holds:

$$\frac{1}{k+1}\left({}^A I_{a^+}^{\alpha,\beta}\theta(\theta+k)(x)\right) \leq \frac{1}{2(k+2)}\left({}^A I_{a^+}^{\alpha,\beta}(2\theta+k)^2\right)(x) \leq \left({}^A I_{a^+}^{\alpha,\beta}\theta(\theta+k)(x)\right)$$

*Proof.* Choosing $f_1(\theta) = \theta + k$ and $f_2(\theta) = k$, we have respectively $\tau_1 = 1$ and $\tau_2 = k+1$

Applying Theorem 4.5, we get the desired result.

**Example 5.8**. Suppose that, $\alpha, \beta, \hbar > 0, p \geq 1$. Let $f_1, f_2 \in L[a, x]$ be two positive functions on $[0, \infty)$ such that for all $x > a \geq 1$ and $\theta \in [a, x]$, the following inequality holds:

$$\left( {}^{A}I_{a^+}^{\alpha,\beta} \theta^p(x) \right)^{\frac{1}{p}} + \left( {}^{A}I_{a^+}^{\alpha,\beta} (\theta + \hbar)^p(x) \right)^{\frac{1}{p}} \leq 2\left( {}^{A}I_{a^+}^{\alpha,\beta} [\Upsilon_\hbar(\theta)]^p \right)^{\frac{1}{p}}(x),$$

where $\Upsilon_\hbar(\theta) = \max\{\hbar(2 + \hbar) + \theta, \theta(1 + \hbar) - \hbar\}$.

*Proof.* Taking $f_1(\theta) = \theta + \hbar$ and $f_2(\theta) = \hbar$, we have respectively $\tau_1 = 1$ and $\tau_2 = \hbar + 1$

Applying Theorem 4.6, we get the desired result.

## 6. Conclusin

In this paper we study some fractional integral inequalities in a generalised sense. Then we obtain inverse Minkowski inequalities for fractional integral operators with a new analytic kernel. Finally, some concrete examples show the importance of our results. Not only do we prove that it is mathematically more valuable than writing individual papers on each operator, but our results also encourage such research for the future. We hope that our results can stimulate further research in various fields of pure and applied science.